\documentclass[a4paper]{article}
\usepackage{amsmath}
\usepackage{amssymb}
\usepackage[square,numbers]{natbib}
\usepackage{url}
\usepackage{csquotes}

\title{$q$-analogs of group divisible designs}
\author{Marco Buratti\thanks{Dipartimento di Matematica e Informatica,
     Universita degli Studi di Perugia, 06123 Perugia, Italy, \texttt{buratti@dmi.unipg.it}},
Michael Kiermaier\thanks{Department of Mathematics, University of Bayreuth, D-95440 Bayreuth, Germany,
	\texttt{michael.kiermaier@uni-bayreuth.de}},
Sascha Kurz\thanks{Department of Mathematics, University of Bayreuth, D-95440 Bayreuth, Germany,
	\texttt{sascha.kurz@uni-bayreuth.de}},
\\
Anamari Naki\'{c}\thanks{Faculty of Electrical Engineering and Computing, University of Zagreb,
	Unska 3, 10000 Zagreb, Croatia,	\texttt{anamari.nakic@fer.hr}},
Alfred Wassermann\thanks{Department of Mathematics, University of Bayreuth, D-95440 Bayreuth, Germany,
	\texttt{alfred.wassermann@uni-bayreuth.de}}}

\newtheorem{definition}{Definition}
\newtheorem{example}{Example}
\newtheorem{theorem}{Theorem}
\newtheorem{lemma}{Lemma}
\newtheorem{remark}{Remark}

\newenvironment{Proof}
{\begin{trivlist}\item[]{{\sc Proof.}}}{\hfill{$\square$}\noindent\end{trivlist}}

\newcommand{\cB}{\mathcal{B}}
\newcommand{\Des}{\mathcal{D}}
\newcommand{\cG}{\mathcal{G}}
\newcommand{\cP}{\mathcal{P}}
\newcommand{\cV}{V}
\newcommand{\delt}{\varDelta}
\newcommand{\qbin}[3]{\genfrac{[}{]}{0pt}{}{#1}{#2}_{#3}}

\DeclareMathOperator{\GF}{GF}
\DeclareMathOperator{\GL}{GL}
\DeclareMathOperator{\SL}{SL}
\DeclareMathOperator{\PG}{PG}
\DeclareMathOperator{\PGGL}{P\Gamma L}

\begin{document}
\maketitle

\begin{abstract}
A well known class of objects in combinatorial design theory are {group divisible designs}.
Here, we introduce the $q$-analogs of group divisible designs.
It turns out that there are interesting connections to
scattered subspaces, $q$-Steiner systems, packing designs and $q^r$-divisible projective sets.

We give necessary conditions for the existence of $q$-analogs of group divisible designs, 
construct an infinite series of examples, and provide further existence results with the help of 
a computer search.

One example is a $(6,2,3,2)_2$ group divisible design over $\GF(2)$ which is a packing design consisting of $180$ blocks that such every $2$-dimensional subspace in $\GF(2)^6$ is covered at most twice.
\end{abstract}
\section{Introduction}

The classical theory of \emph{$q$-analogs} of mathematical objects and functions has
its beginnings as early as in the work of Euler~\cite{Euler1753}. 
In 1957, Tits~\cite{Tits-1957} further suggested
that combinatorics of sets could be regarded as the limiting case $q \to 1$ of combinatorics 
of vector spaces over the finite field $\GF(q)$. 
Recently, there has been an increased interest in studying $q$-analogs of combinatorial designs 
from an applications' view. These $q$-analog structures can be useful in 
network coding and distributed storage, see e.g.~\cite{networkcoding2018}. 

It is therefore natural to ask
which combinatorial structures can be generalized from sets to
vector spaces over $\GF(q)$. 
For combinatorial designs, this question was first studied by 
Ray-Chaudhuri ~\cite{Berge1974}, Cameron~\cite{Cameron-1974,Cameron2} and Delsarte~\cite{Delsarte-1976-JCTSA20[2]:230-243} in the early 1970s. 

Specifically, let $\GF(q)^v$ be a vector space of dimension $v$ over the finite field $\GF(q)$. 
Then a $t\text{-}(v, k, \lambda)_q$ subspace design is defined as a collection of $k$-dimensional subspaces
of $\GF(q)^v$, called blocks, such that each $t$-dimensional subspace of $\GF(q)^v$ is contained
in exactly $\lambda$ blocks. Such $t$-designs over $\GF(q)$ are the $q$-analogs of conventional
designs. By analogy with the $q \to 1$ case, a $t\text{-}(v, k, 1)_q$ subspace design is said to
be a $q$-Steiner system, and denoted by $S(t, k, v)_q$.

Another well-known class of objects in combinatorial design theory are \emph{group divisible designs}~\cite{Mullin-Gronau-2007}. 
Considering the above, it therefore seems natural to ask for $q$-analogs of group divisible designs.

Quite surprisingly, it turns out that $q$-analogs of group divisible designs
have interesting connections to scattered subspaces which are central objects in finite geometry,
as well as to coding theory via  $q^r$-divisible projective sets.
We will also discuss the connection to $q$-Steiner systems~\cite{fmp:10491987} and to 
packing designs~\cite{EtzionZhang2018}. 

\medskip
Let $k$, $g$, and $\lambda$ be positive integers.
A $(v, g, k, \lambda)$-\emph{group divisible design}
of index $\lambda$ and order $v$ is a triple $(\cV, \cG, \cB)$, where
$\cV$ is a finite set of
cardinality $v$, $\cG$, where $\#  \cG > 1$,  is a partition of $\cV$ into parts (groups) of cardinality $g$, and
$\cB$ is a family of subsets (blocks) of $\cV$ (with $\# B = k$ for $B\in \cB$)
such that every pair of distinct elements of $\cV$ occurs in exactly  $\lambda$ blocks or one group, but not both.

See---for example---\cite{HANANI1975255,Mullin-Gronau-2007} for details.
We note that the ``groups'' in group divisible designs have nothing to do with group theory.

The $q$-analog of a combinatorial structure over sets
is defined by replacing subsets by subspaces and cardinalities by dimensions.
Thus, the $q$-analog of a group divisible design can be defined as follows.
\begin{definition}\label{def:qgdd}
Let $k$, $g$, and $\lambda$ be positive integers.
A \emph{$q$-analog of a group divisible design} of index $\lambda$ and
order $v$ --- denoted as \emph{$(v, g, k, \lambda)_q$-GDD} ---
is a triple $(\cV, \cG, \cB)$, where
\begin{itemize}
\item[--] $\cV$ is a vector space over $\GF(q)$ of dimension $v$,
\item[--] $\cG$ is a vector space partition\footnote{A set of subspaces of $\cV$ such that every $1$-dimensional subspace is 
covered exactly once is called vector space partition.} of $\cV$ into subspaces (groups) of dimension $g$, and
\item[--] $\cB$ is a family  of subspaces (blocks) of $\cV$,
\end{itemize}
that satisfies
\begin{enumerate}
\item $\#  \cG > 1$,
\item if $B\in \cB$ then $\dim B = k$,
\item every $2$-dimensional subspace of $\cV$ occurs in exactly  $\lambda$ blocks or one group, but not both.
\end{enumerate}
\end{definition}

In the sequel, we will only consider so called \emph{simple} group divisible designs, 
i.e.~designs without multiple appearances of blocks.

In finite geometry a partition of the $1$-dimensional subspaces of $\cV$ in subspaces of dimension $g$ 
is known as \emph{$(g-1)$-spread}.

This notation respects the well-established usage of the geometric dimension $(g-1)$ of the spread elements.
Nevertheless, for the rest of the paper we think of the elements of a $(g-1)$-spread 
as subspaces of algebraic dimension $g$ of a $v$-dimensional vector space $V$. 
Similarly, $2$-dimensional subspaces of $V$ will sometimes be called \emph{lines}.

\medskip
A possible generalization would be to require the last condition in Definition~\ref{def:qgdd}  
for every $t$-dimensional subspace of $\cV$, where $t\ge 2$. 
For $t=1$ such a definition would make no sense. 

An equivalent formulation of the last condition in Definition~\ref{def:qgdd} 
would be that every block in $\cB$ intersects the spread elements in dimension of at most one.
The $q$-analog of concept of a \emph{transversal design} would be that 
every block in $\cB$ intersects the spread elements exactly in dimension one.
But for $q$-analogs this is only possible in the trivial case $g=1$, $k=v$. 
However, a related concept was defined in \cite{etzion2013codes}.

\smallskip
Another generalization of Definition~\ref{def:qgdd} which is well known for the set case is:

Let $K$ and $G$ be sets of positive integers and let $\lambda$ be a positive integer.
A triple $(\cV, \cG, \cB)$ is called \emph{$(v, G, K, \lambda)_q$-GDD}, if
$\cV$ is a vector space over $\GF(q)$ of dimension $v$,
$\cG$ is a vector space partition of $\cV$ into subspaces (groups) whose dimensions lie in $G$, and
$\cB$ is a family  of subspaces (blocks) of $\cV$,
that satisfies
\begin{enumerate}
\item $\#  \cG > 1$,
\item if $B\in \cB$ then $\dim B \in K$,
\item every $2$-dimensional subspace of $\cV$ occurs in exactly $\lambda$ blocks or one group, but not both.
\end{enumerate}
Then, a $(v, \{g\}, K, \lambda)_q$-GDD is called $g$-\emph{uniform}.

\smallskip
An even more general definition --- which is also studied in the set case --- is 
a \emph{$(v, G, K, \lambda_1, \lambda_2)_q$-GDD} for which condition 3.~is replaced by
\begin{enumerate}
\item[3'.] every $2$-dimensional subspace of $\cV$ occurs in $\lambda_1$ blocks if it is contained in a group, otherwise
it is contained in exactly $\lambda_2$ blocks.
\end{enumerate}
Thus, a $q$-GDD of Definition~\ref{def:qgdd} is a $(v, \{g\}, \{k\}, 0, \lambda)_q$-GDD in the general form.

\medskip
Among all $2$-subspaces of $V$, only a small fraction is covered by the elements of $\cG$.
Thus, a $(v, g, k, \lambda)_q$-GDD is ``almost'' a $2$-$(v,k,\lambda)_q$ subspace design, 
in the sense that the vast majority of the $2$-subspaces is covered by 
$\lambda$ elements of $\cB$.
From a slightly different point of view, a $(v,g,k,\lambda)_q$-GDD is a 
$2$-$(v,g,k,\lambda)_q$ \emph{packing design} of fairly large size, 
which are designs where the condition 
``each $t$-subspace is covered by exactly $\lambda$ blocks'' is relaxed to 
``each $t$-subspace is covered by at most $\lambda$ blocks'' \cite{qdesignscomputer2017}.
In Section~\ref{sec:computer} we give an example of a 
$(6,2,3,2)_2$-GDD consisting of $180$ blocks. 
This is the largest known $2$-$(6,3,2)_2$ packing design.

We note that a $q$-analog of a group divisible design
can be also seen as a special graph decomposition over
a finite field, a concept recently introduced in~\cite{BNW:2018}.
It is indeed equivalent to a decomposition of a complete $m$-partite graph into cliques where:
the vertices are the points of a projective space $\PG(n,q)$;
the parts are the members of a spread of $\PG(n,q)$ into subspaces of a suitable dimension;
the vertex-set of each clique is a subspace of $\PG(n,q)$ of a suitable dimension.

\section{Preliminaries}
For $1\leq m \leq v$ we denote the set of $m$-dimensional subspaces of $\cV$, also called
\emph{Grassmannian}, by $\qbin{\cV}{m}{q}$. It is well known that its cardinality can be expressed by the
\emph{Gaussian coefficient}
\[
	\#\qbin{\cV}{m}{q} = \qbin{v}{m}{q} = \frac{(q^v-1)(q^{v-1}-1)\cdots(q^{v-m+1}-1)}{(q^m-1)(q^{m-1}-1)\cdots(q-1)}\,.
\]

\begin{definition}
Given a spread in dimension $v$, let $\qbin{\cV}{k}{q}'$
be the set of all $k$-dimensional subspaces in $\cV$ that
contain no $2$-dimensional subspace which is already covered by the spread.
\end{definition}
The intersection between a $k$-dimensional subspace $B\in\qbin{\cV}{2}{q}'$
and each element of the spread is at most one-dimensional.
In finite geometry such a subspace $B\in\qbin{\cV}{k}{q}'$ is called
\emph{scattered subspace with respect to $\cG$}, see~\cite{BALL2000294,Blokhuis2000}.

In case $g=1$, i.e. $\cG = \qbin{\cV}{1}{q}$, no $2$-dimensional subspace is
covered by this trivial spread. Then, $(\cV, \cB)$ is a 
$2$-$(v,k,\lambda)_q$ subspace design. 
See \cite{qdesignscomputer2017,qdesigns2017} for surveys about subspace designs 
and computer methods for their construction.


\medskip
Let $g\cdot s = v$ and $\cV = \GF(q)^{v}$.
Then, the set of $1$-dimensional subspaces of $\GF(q^g)^s$ regarded as $g$-dimensional subspaces
in the $q$-linear vector space $\GF(q)^{v}$, i.e.
\[\cG = \qbin{\GF(q^g)^s}{1}{q^g}\,,
\]
is called \emph{Desarguesian spread}.

A $t$-spread $\cG$ is called \emph{normal} or \emph{geometric}, 
if $U,V\in \cG$ then any element $W\in\cG$ is either
disjoint to the subspace $\langle U,V\rangle$ or contained in it, see e.g.~\cite{Lunardon1999}.
Since all normal spreads are isomorphic to the Desarguesian spread~\cite{Lunardon1999}, we will 
follow \cite{Lavrauw2016} and denote normal spreads as Desarguesian spreads.

If $s\in\{1,2\}$, then all spreads are normal and therefore Desarguesian. 
The automorphism group of a Desarguesian spread $\cG$ is $\PGGL(s, q^g)$. 


\paragraph{``Trivial'' $q$-analogs of group divisible designs.} 
For subspace designs, the empty set as well as the the set of all 
$k$-dimensional subspaces in $\GF(q)^v$ always
are designs, called \emph{trivial designs}.
Here, it turns out that the question if trivial $q$-analogs of group divisible designs exist
is rather non-trivial.

Of course, iff $g\mid v$, there exists always the trivial $(v,g,k,0)_q$-GDD $(V,\cG, \{\})$.
But it is not clear if the set of all scattered $k$-dimensional subspaces, 
i.e. $(V,\cG, \qbin{\cV}{k}{q}')$, is always a $q$-GDD.
This would require that 
every subspace $L\in\qbin{\cV}{2}{q}$ that is not covered by the spread,
is contained in the same number $\lambda_{\max}$ of blocks of $\qbin{\cV}{k}{q}'$. 
If this is the case, we call $(\cV, \qbin{\cV}{k}{q}', \cG)$ the \emph{complete} 
$(v, g, k, \lambda_{\max})_q$-GDD.

If the complete $(v, g, k, \lambda_{\max})_q$-GDD exists, 
then for any $(v, g, k, \lambda)_q$-GDD $(\cV, \cG, \cB)$ the triple
$(\cV, \cG, \qbin{\cV}{2}{q}'\setminus\cB)$ is a $(v, g, k, \lambda_{\max}-\lambda)_q$-GDD,
called the \emph{supplementary} $q$-GDD.


For a few cases we can answer the question if the complete $q$-GDD exists, or in other words, if there 
is a $\lambda_{\max}$.
In general, the answer depends on the choice of the spread.
In the smallest case, $k=3$, however, $\lambda_{\max}$ exists for all spreads.
\begin{lemma}
Let $\cG$ be a $(g-1)$-spread in $V$ and let
$L$ be a $2$-dimensional subspace which is not contained in any element of $\cG$.
Then, $L$ is contained
in
\[
	\lambda_{\max} = \qbin{v-2}{3-2}{q} -
		\qbin{2}{1}{q} \qbin{g-1}{3-2}{q}
\]
blocks of $\qbin{\cV}{3}{q}'$.
\end{lemma}
\begin{Proof}
Every $2$-dimensional subspace $L$ is contained in
$\qbin{v-2}{3-2}{q}$ $3$-dimensional subspaces of $V$.
If $L$ is not contained in any spread element, 
this means that $L$ intersects $\qbin{2}{1}{q}$ different spread elements and the intersections
are $1$-dimensional.
Let $S$ be one such spread element.
Now, there are $\qbin{g-1}{1}{q}$ choices among the $3$-dimensional subspaces
in $\qbin{V}{3}{q}$ which contain $L$ to
intersect $S$ in dimension two.
Therefore, $L$ is contained in 
\[
	\lambda_{\max} = \qbin{v-2}{3-2}{q} -
		\qbin{2}{1}{q} \qbin{g-1}{3-2}{q}
\]
blocks of $\qbin{\cV}{3}{q}'$.
\end{Proof}

In general, the existence of $\lambda_{\max}$ may depend on the spread. 
This can be seen 
from the fact that the maximum dimension of a scattered subspace 
depends on the spread, see \cite{Blokhuis2000}.
However, for a Desarguesian spread and $g=2$, $k=4$, we can determine $\lambda_{\max}$.
\begin{lemma}
Let $\cG$ be a Desarguesian $1$-spread in $V$ and let
$L$ be a $2$-dimensional subspace which is not contained in any element of $\cG$.
Then, $L$ is contained
in
\[
	\lambda_{\max} = \qbin{v-2}{4-2}{q} -
		1-q\qbin{2}{1}{q}\qbin{v-4}{1}{q} 
		- \qbin{v}{1}{q}/\qbin{2}{1}{q}+\qbin{4}{1}{q}/\qbin{2}{1}{q}
\]
blocks of $\qbin{\cV}{4}{q}'$.
\end{lemma}
\begin{Proof}
Every $2$-dimensional subspace $L$ is contained in
$\qbin{v-2}{4-2}{q}$ $4$-dimensional subspaces.
If $L$ is not covered by the spread
this means that $L$ intersects $\qbin{2}{1}{q}$ spread elements $S_1,\dots,S_{q+1}$, 
which span a subspace $F$. Since the spread is Desarguesian, the dimension of $F$ is equal to $4$.
All other spread elements are disjoint to $L$. Since $L\le F$, 
we have to subtract one possibility. 
For each $1\le i\le q+1$, $\langle S_i,L\rangle$ is contained 
in $q\qbin{v-4}{1}{q}$ $4$-dimensional subspaces with a $3$-dimensional intersection with $F$. 
All other spread elements $S'$ of $F$ satisfy $\langle S',L\rangle=F$. 
If $S''$ is one of the 
$\qbin{v}{1}{q}/\qbin{2}{1}{q}-\qbin{4}{1}{q}/\qbin{2}{1}{q}$  spread elements disjoint from $F$, then 
$F'':=\langle S'',L\rangle$ intersects $F$ in dimension $2$. Moreover, $F''$ does not contain any further 
spread element, since otherwise $F''$ would be partitioned into $q^2+1$ spread elements, where $q+1$ of them 
have to intersect $L$. Thus, $L$ is contained in exactly $\lambda_{\max}$ elements 
from $\qbin{\cV}{4}{q}'$.
\end{Proof}

\section{Necessary conditions on $(v, g, k,\lambda)_q$}
The necessary conditions for a $(v, g, k, \lambda)$-GDD over sets are
    $g \mid v$,
    $k \leq v/g$,
    $\lambda (\frac{v}{g} - 1) g \equiv 0 \pmod{k-1}$,
    and
    $\lambda \frac{v}{g}(\frac{v}{g} - 1) g^2 \equiv 0 \pmod{k(k-1)}$,
see \cite{HANANI1975255}.

For $q$-analogs of GDDs it is well known that $(g-1)$-spreads exist if and only if $g$ divides $v$.
A $(g-1)$-spread consists of $\qbin{v}{1}{q} / \qbin{g}{1}{q}$ blocks and contains
\[
\qbin{g}{2}{q}\cdot \qbin{v}{1}{q} / \qbin{g}{1}{q}
\]
$2$-dimensional subspaces.

Based on the pigeonhole principle we can argue that if $B$ is a block of a $(v, g, k, \lambda)_q$ $q$-GDD 
then there cannot be more points in $B$ than the number of spread elements,
i.e.~if
$\qbin{k}{1}{q} \leq \qbin{v}{1}{q} / \qbin{g}{1}{q}$.
It follows that (see \cite[Theorem 3.1]{Blokhuis2000})
\begin{equation}\label{eq:pigeon}
k\leq v - g\,.
\end{equation}
This is the $q$-analog of the restriction $k\leq v/g$ for the set case.

If $\cG$ is a Desarguesian spread, it follows
from \cite[Theorem 4.3]{Blokhuis2000}
for the parameters $(v, g, k, \lambda)_q$ to be admissible
that
\[
	k \leq v / 2\,.
\]

By looking at the numbers of $2$-dimensional subspaces which are covered by spread elements
we can conclude that the cardinality of $\cB$ has to be
\begin{equation}\label{eq:nec}
\#\cB = \lambda\frac{\qbin{v}{2}{q} - \qbin{g}{2}{q}\cdot \qbin{v}{1}{q} / \qbin{g}{1}{q}}{\qbin{k}{2}{q}} \,.
\end{equation}
A necessary condition on the parameters of a $(v, g, k, \lambda)_q$-GDD is that
the cardinality in (\ref{eq:nec}) is an integer number.

Any fixed $1$-dimensional subspace $P$ is contained in $\qbin{v-1}{1}{q}$ $2$-dimensional subspaces.
Further, $P$ lies in exactly one block of the spread and this block covers
$\qbin{g-1}{1}{q}$ $2$-dimensional subspaces through $P$.
Those $2$-dimensional subspaces are not covered by blocks in $\cB$.
All other $2$-dimensional subspaces containing $P$ are covered by exactly $\lambda$ $k$-dimensional
blocks. Such a block contains $P$ and there are $\qbin{k-1}{1}{q}$ $2$-dimensional subspaces through $P$
in this block.
It follows that $P$ is contained in exactly
\begin{equation}\label{eq:nec2}
\lambda\frac{\qbin{v-1}{1}{q} - \qbin{g-1}{1}{q}}{\qbin{k-1}{1}{q}}
\end{equation}
$k$-dimensional blocks and this number must be an integer.
The number (\ref{eq:nec2}) is the \emph{replication number}
of the point $P$ in the $q$-GDD.

\medskip
Up to now, the restrictions
(\ref{eq:pigeon}), (\ref{eq:nec}), (\ref{eq:nec2}), as well as $g$ divides $v$,
on the parameters of a $(v, g, k,\lambda)_q$-GDD
are the $q$-analogs of restrictions for the set case.
But for $q$-GDDs there is a further necessary condition whose analog
in the set case is trivial.


Given a multiset of subspaces of $\cV$,
we obtain a corresponding multiset $\cP$ of points by
replacing each subspace by its set of points.
A multiset $\cP\subseteq\qbin{\cV}{1}{q}$ of points
in $\cV$ can be expressed by its weight function $w_\cP$:
For each point $P\in\cV$ we denote its multiplicity in $\cP$ by $w_\cP(P)$.
We write
\[
\# \cP = \sum_{P\in \cV} w_\cP(P)
\quad\mbox{ and  }\quad
\# (\cP \cap H)=\sum_{P\in H} w_\cP(P)
\]
where $H$ is an arbitrary hyperplane in $\cV$.

Let
$1\le r<v$ be an integer.
If
$\#\cP\equiv \# (\cP\cap H) \pmod {q^{r}}$ for every hyperplane $H$,
then $\cP$ is called \emph{$q^r$-divisible}.\footnote{Taking the elements of $\cP$ as columns of a generator matrix
gives a linear code of length $\#\cP$ and dimension $k$ whose codewords have weights being divisible by $q^r$.}
In \cite[Lemma 1]{upper_bounds_cdc} it is shown that the
multiset $\cP$ of points corresponding to a multiset of subspaces with dimension at least $k$ is $q^{k-1}$-divisible.
\begin{lemma}[{\cite[Lemma 1]{upper_bounds_cdc}}]
For a non-empty multiset of subspaces of $\cV$ with $m_i$ subspaces of dimension $i$ let $\cP$ be the corresponding
multiset of points.
If $m_i=0$ for all $0\leq i < k$, where $k\geq 2$, then
\[
\# \cP \equiv \#(\cP\cap H) \pmod{q^{k-1}}
\]
for every hyperplane $H\leq\cV$.
\end{lemma}
\begin{Proof}
We have $\#\cP = \sum_{i=0}^vm_i\qbin{v}{i}{q}$.
The intersection of an $i$-subspace $U \leq \cV$ with an arbitrary hyperplane $H\leq\cV$
has either dimension $i$ or $i-1$.
Therefore, for the set $\cP'$ of points corresponding to $U$, we get that
$\#\cP = \qbin{i}{1}{q}$ and that $\#(\cP'\cap H)$ is equal to $\qbin{i}{1}{q}$ or $\qbin{i-1}{1}{q}$.
In either case, it follows from
$\qbin{i}{1}{q}\equiv \qbin{i-1}{1}{q} \pmod{q^{i-1}}$
that
\[
\#(\cP'\cap H) \equiv \qbin{i}{1}{q} \pmod{q^{i-1}} \,.
\]
Summing up yields the proposed result.
\end{Proof}
If there is a suitable integer $\lambda$ such that
$w_{\cP}(P)\le \lambda$ for all $P\in\cV$ , then we can define for $\cP$ the complementary weight function
\[
\bar{w}_\lambda(P) = \lambda - w(P)
\]
which in turn gives rise to the \emph{complementary} multiset of points $\bar{\cP}$.
In \cite[Lemma 2]{upper_bounds_cdc} it is shown that a $q^r$-divisible multiset $\cP$ leads to a multiset
$\overline{\cP}$ that is also $q^r$-divisible.
\begin{lemma}[{\cite[Lemma 2]{upper_bounds_cdc}}]
If a multiset $\cP$ in $\cV$ is $q^r$-divisible with $r<v$ and
satisfies $w_{\cP}(P)\le \lambda$ for all $P\in\cV$
then the complementary multiset $\bar{\cP}$
is also $q^r$-divisible.
\end{lemma}
\begin{Proof}
We have
\[
\#\bar{\cP}  = \qbin{v}{1}{q} \lambda - \#\cP
\quad\mbox{ and }\quad
\#(\bar{\cP}\cap H)  = \qbin{v-1}{1}{q} \lambda - \#(\cP\cap H)
\]
for every hyperplane $H\leq \cV$.
Thus, the result follows from $\qbin{v}{1}{q}\equiv \qbin{v-1}{1}{q} \pmod{q^r}$
which holds for $r<v$.
\end{Proof}
These easy but rather generally applicable facts about
$q^r$-divisible multiset of points are enough to conclude:
\begin{lemma}\label{lemma_div}
  Let $(\cV, \cG, \cB)$ be a $(v, g, k,\lambda)_q$-GDD and
  $2\leq g \leq k$,
  then $q^{k-g}$ divides $\lambda$.
\end{lemma}
\begin{Proof}
  Let $P\in \qbin{\cV}{1}{q}$ be an arbitrary point.
  Then there exists exactly one spread element $S\in \cG$ that contains $P$.
  By $\cB_P$ we denote the elements of $\cB$ that contain $P$.
  Let $S'$ and $\cB_P'$ denote the corresponding subspaces in the factor space $\cV / P$.

  We observe that every point of $\qbin{S'}{1}{q}$ is disjoint to the elements of $\cB_P'$ and that
  every point in $\qbin{\cV / P}{1}{q} \setminus \qbin{S'}{1}{q}$ is met by exactly $\lambda$ elements of $\cB_P'$
  (all having dimension $k-1$).
  We note that $\cB_P'$ gives rise to a $q^{k-2}$-divisible multiset $\cP$ of points.
  So, its complement $\bar{\cP}$, which is the $\lambda$-fold copy of $S'$,
  also has to be $q^{k-2}$-divisible.
  For every hyperplane $H$ not containing $S'$, we have $\#(\bar{\cP}\cap H)=\lambda\qbin{g-2}{1}{q}$ and 
  $\#\bar{\cP}=\lambda\qbin{g-1}{1}{q}$. Thus, $\lambda q^{g-2}=\#\bar{\cP}-\#(\bar{\cP}\cap H)\equiv 0\pmod{q^{k-2}}$, 
  so that $q^{k-g}$ divides $\lambda$.  
\end{Proof}

We remark that the criterion in Lemma~\ref{lemma_div} is independent of the dimension $v$ of the ambient space.
Summarizing the above we arrive at the following restrictions.
\begin{theorem}
Necessary conditions for a $(v, g, k, \lambda)_q$-GDD are
\begin{enumerate}
\item $g$ divides $v$,
\item $k\leq v - g$,
\item the cardinalities in
(\ref{eq:nec}), (\ref{eq:nec2}) are integer numbers,
\item if $2 \leq g \leq k$ then $q^{k-g}$ divides $\lambda$.
\end{enumerate}
If these conditions are fulfilled,
the parameters $(v, g, k, \lambda)_q$ are called \emph{admissible}.
\end{theorem}
 Table~\ref{tab:nec} contains the admissible parameters for $q=2$ up to dimension $v=14$.
Column $\lambda_\delt$ gives the minimum value of $\lambda$
which fulfills the above necessary conditions.
All admissible values of $\lambda$ are integer multiples of
$\lambda_\delt$.
In column $\# \cB$ the cardinality of $\cB$ is given for 
$\lambda=\lambda_\delt$. Those values of $\lambda_{\max}$ that are valid for the  
Desarguesian spread only are given in italics, where the values for 
$(v,g,k)=(8,4,4)$ and $(9,3,4)$ have been checked by a computer enumeration.


\begin{table}
\caption{Admissible  parameters for $(v, g, k,\lambda)_2$-GDDs with $v\leq 14$.}\label{tab:nec}
\begin{center}
{\small
\begin{tabular}{rrrrrrr}
$v$ & $g$ & $k$ & $\lambda_{\delt}$ & $\lambda_{\max}$ & $\#\cB$ & $\# \cG$ \\  \hline
6 & 2 & 3 & 2 & 12 & 180 & 21\\
6 & 3 & 3 & 3 & 6 & 252 & 9\\
8 & 2 & 3 & 2 & 60 & 3060 & 85\\
8 & 2 & 4 & 4 & \textit{480} & 1224 & 85\\
8 & 4 & 3 & 7 & 42 & 10200 & 17\\
8 & 4 & 4 & 7 & \textit{14}  & 2040 & 17\\
9 & 3 & 3 & 1 & 118 & 6132 & 73\\
9 & 3 & 4 & 10 & \textit{1680}   & 12264 & 73\\
10 & 2 & 3 & 14 & 252 & 347820 & 341\\
10 & 2 & 4 & 28 & \textit{10080} & 139128 & 341\\
10 & 2 & 5 & 8 &    & 8976 & 341\\
10 & 5 & 3 & 21 & 210 & 507408 & 33\\
10 & 5 & 4 & 35 &    & 169136 & 33\\
10 & 5 & 5 & 15 &    & 16368 & 33\\
12 & 2 & 3 & 2 & 1020 & 797940 & 1365\\
12 & 2 & 4 & 28 & \textit{171360}  & 2234232 & 1365\\ 
12 & 2 & 5 & 40 &    & 720720 & 1365\\
12 & 2 & 6 & 16 &    & 68640 & 1365\\
12 & 3 & 3 & 3 & 1014 & 1195740 & 585\\
12 & 3 & 4 & 2 &    & 159432 & 585\\
12 & 3 & 5 & 1860 &    & 33480720 & 585\\
12 & 3 & 6 & 248 &    & 1062880 & 585\\
12 & 4 & 3 & 1 & 1002 & 397800 & 273\\
12 & 4 & 4 & 7 &    & 556920 & 273\\
12 & 4 & 5 & 62 &    & 1113840 & 273\\
12 & 4 & 6 & 124 &    & 530400 & 273\\
12 & 6 & 3 & 1 & 930 & 393120 & 65\\
12 & 6 & 4 & 1 &    & 78624 & 65\\
12 & 6 & 5 & 155 &    & 2751840 & 65\\
12 & 6 & 6 & 31 &    & 131040 & 65\\
14 & 2 & 3 & 2 & 4092 & 12778740 & 5461\\
14 & 2 & 4 & 4 & \textit{2782560} & 5111496 & 5461\\ 
14 & 2 & 5 & 248 &    & 71560944 & 5461\\
14 & 2 & 6 & 496 &    & 34076640 & 5461\\
14 & 2 & 7 & 32 &    & 536640 & 5461\\
14 & 7 & 3 & 21 & 3906 & 133161024 & 129\\
14 & 7 & 4 & 35 &    & 44387008 & 129\\
14 & 7 & 5 & 465 &    & 133161024 & 129\\
14 & 7 & 6 & 651 &    & 44387008 & 129\\
14 & 7 & 7 & 63 &    & 1048512 & 129\\
\end{tabular}}
\end{center}
\end{table}

For the case $\lambda=1$, the online tables \cite{ubt_epub2670}
\begin{center}
\texttt{http://subspacecodes.uni-bayreuth.de}
\end{center}
may give further restrictions, since $\cB$ is
a constant dimension subspace code of minimum distance $2(k-1)$ and therefore
\[
\#\cB \leq A_q(v, 2(k-1); k).
\]
The currently best known upper bounds for $A_q(v,d;k)$ are given by
\cite[Equation~(2)]{a_2_8_6_4_257} referring back to partial spreads and $A_2(6,4;3)=77$ \cite{honold2015optimal}, 
$A_2(8,6;4)=257$ \cite{a_2_8_6_4_257} both obtained by exhaustive integer linear programming computations,
see also \cite{upper_bounds_cdc}.

\section{$q$-GDDs and $q$-Steiner systems}
In the set case the connection between Steiner systems $2$-$(v,k,1)$ and group divisible designs is well understood.
\begin{theorem}[{\cite[Lemma 2.12]{HANANI1975255}}]\label{thm:hanani}
A $2$-$(v+1, k, 1)$ design exists if and only if a $(v, k-1, k, 1)$-GDD exists.
\end{theorem}

There is a partial $q$-analog of Theorem~\ref{thm:hanani}:
\begin{theorem}
\label{thm:qhanani}
If there exists a $2$-$(v+1, k, 1)_q$ subspace design,
then a $(v, k-1, k, q^2)_q$-GDD exists.
\end{theorem}

\begin{Proof}
Let $\cV'$ be a vector space of dimension $v+1$ over $\GF(q)$.
We fix a point $P\in\qbin{\cV'}{1}{q}$ and define the projection
\[
\pi : \operatorname{PG}(\cV') \to \operatorname{PG}(\cV'/P), \quad U \mapsto (U+P)/P\text{.}
\] 
\noindent
For any subspace $U\leq \cV'$ we have
\[
	\dim(\pi(U))
	= \begin{cases}
		\dim(U) - 1 & \text{if }P \leq U\text{,} \\
		\dim(U) & \text{otherwise.}
	\end{cases}
\]

Let $\Des = (\cV', \cB')$ be a $2$-$(v+1, k, 1)_q$ subspace design.
The set
\[
	\cG = \{\pi(B) \mid B\in \cB', P\in B \}
\]
is the derived design of $\Des$ with respect to $P$ \cite{Kiermaier-Laue-2015-AiMoC9[1]:105-115}, which has the parameters $1$-$(v, k-1, 1)_q$.
In other words, it is a $(k-2)$-spread in $\cV'/P$.
Now define
\[
	\cB = \{\pi(B) \mid B\in \cB', P\notin B \} \mbox{ and }  \cV = \cV'/P\,.
\]
We claim that $(\cV, \cG, \cB)$ is a $(v, k-1, k, q^2)_q$-GDD.


In order to prove this, let $L\in\qbin{V}{2}{q}$ be a line not covered by any element in $\cG$.
Then $L = E/P$, where $E\in\qbin{\cV'}{3}{q}$, $P\le E$ and $E$ is not contained in a block of the design $\Des$.
The blocks of $\cB$ covering $L$ have the form $\pi(B)$ with $B\in\cB'$ such that $B \cap E$ is  a line in $E$ not passing through $P$.
There are $q^2$ such lines and each line is contained in a unique block in $\cB'$.
Since these $q^2$ blocks $B$ have to be pairwise distinct and do not contain the point $P$, we get that there are $q^2$ blocks $\pi(B)\in\cB$ containing $L$.
\end{Proof}

Since there are $2$-$(13, 3, 1)_2$ subspace designs~\cite{fmp:10491987}, by Theorem~\ref{thm:qhanani} there are also $(12, 2, 3, 4)_2$-GDDs.

The smallest admissible case of a $2$-$(v,3,1)_q$ subspace design is $v=7$, which is known as a \emph{$q$-analog of the Fano plane}.
Its existence is a notorious open question for any value of $q$.
By Theorem~\ref{thm:qhanani}, the existence would imply the existence of a 
$(6, 2, 3, q^2)_q$-GDD, which has been shown to be true in \cite{EtzionHooker2017} 
for any value of $q$, in the terminology of a 
``residual construction for the $q$-Fano plane''.
In Theorem~\ref{thm_qgdd_orbit_construction}, we will give a general construction of $q$-GDDs covering these parameters.
The crucial question is if a $(6, 2, 3, q^2)_q$-GDD can be ``lifted'' to a $2$-$(7, 3, 1)_q$ 
subspace design. 
While the GDDs with these parameters constructed in Theorem~\ref{thm_qgdd_orbit_construction} 
have a large automorphism group, for the binary case $q=2$ we know from 
\cite{Braun-Kiermaier-Nakic-2016-EuJC51:443-457, kiermaier2016order} that the order of the 
automorphism group of a putative $2$-$(7, 3, 1)_2$ subspace design is at most two.
So if the lifting construction is at all possible for the binary $(6,2,3,4)_2$-GDD from
Theorem~\ref{thm_qgdd_orbit_construction}, necessarily many automorphisms have to 
``get destroyed''.

In Table~\ref{tab:results_q2} we can see that there exists a $(8, 2, 3, 4)_2$-GDD.
This might lead in the same way to a $2$-$(9, 3, 1)_2$ subspace design, which is not known to exist.

\section{A general construction}

A very successful approach to construct $t$-$(v,k,\lambda)$ designs over sets is to prescribe an automorphism group which acts transitively on the subsets of cardinality $t$. 
However for $q$-analogs of designs with $t\geq 2$ this approach yields only trivial designs, since in \cite[Prop. 8.4]{cameron1979384} it is shown that if a group $G\leq \PGGL(v,q)$ acts transitively on the $t$-dimensional subspaces of $\cV$, $2\leq t \leq v-2$, then $G$ acts transitively also on the $k$-dimensional subspaces of $\cV$ for all $1\leq k \leq v-1$.

The following lemma provides the counterpart of the construction idea for $q$-analogs of group divisible designs.
Unlike the situation of $q$-analogs of designs, in this slightly different setting there are indeed suitable groups admitting the general construction of non-trivial $q$-GDDs, which will be described in the sequel.
Itoh's construction of infinite families of subspace designs is based on a similar idea~\cite{Itoh-1998}.

\begin{lemma}
	\label{lem:transitive}
	Let $\cG$ be a $(g-1)$-spread in $\PG(V)$ and let $G$ be a subgroup
of the stabilizer $\PGGL(v,q)_{\cG}$ of $\cG$ in $\PGGL(v,q)$.
	If the action of $G$ on $\qbin{V}{2}{q} \setminus \bigcup_{S\in\cG}\qbin{S}{2}{q}$ is transitive, then any union $\cB$ of $G$-orbits on the set of $k$-subspaces which are scattered with respect to $\cG$ yields a $(v,g,k,\lambda)_q$-GDD $(V,\cG,\cB)$ for a suitable value~$\lambda$.
\end{lemma}

\begin{Proof}
	By transitivity, the number $\lambda$ of blocks in $\cB$ passing through a line $L\in\qbin{V}{2}{q} \setminus \bigcup_{S\in\cG}\qbin{S}{2}{q}$ does not depend on the choice of $L$.
\end{Proof}

In the following, let $V = \GF(q^g)^s$, which is a vector space over $\GF(q)$ of dimension $v = gs$.
Furthermore, let $\cG = \qbin{V}{1}{q^g}$ be the Desarguesian $(g-1)$-spread in $\PG(V)$.
For every $\GF(q)$-subspace $U \leq V$ we have that
\[
\dim_{\GF(q^g)}\bigl(\langle U\rangle_{\GF(q^g)}\bigr) \leq \dim_{\GF(q)}(U)\text{.}
\]
In the case of equality, $U$ will be called \emph{fat}.
Equivalently, $U$ is fat if and only if one (and then any) $\GF(q)$-basis of $U$ is $\GF(q^g)$-linearly independent.
The set of fat $k$-subspaces of $V$ will be denoted by $\mathcal{F}_k$.

We remark that for a fat subspace $U$, the set of points $\{\langle x\rangle_{\GF(q^g)} : x \in U\}$ 
is a Baer subspace of $V$ as a $\GF(q^g)$-vector space.

\begin{lemma}
\label{lem:fatsize}
\[
	\#\mathcal{F}_k = q^{(g-1)\binom{k}{2}} \prod_{i=0}^{k-1} \frac{q^{g(s-i)} - 1}{q^{k-i} - 1}\text{.}
\]
\end{lemma}

\begin{Proof}
	A sequence of $k$ vectors in $V$ is the $\GF(q)$-basis of a fat $k$-subspace if and only if it is linearly independent over $\GF(q^g)$.
	Counting the set of those sequences in two ways yields
	\[
		\#\mathcal{F}_k \cdot\prod_{i=0}^{k-1} (q^k - q^i) = \prod_{i=0}^{k-1} ((q^g)^s - (q^g)^i)\text{,}
	\]
	which leads to the stated formula.
\end{Proof}

We will identify the unit group $\GF(q)^*$ with the corresponding group of $s\times s$ scalar matrices over $\GF(q^g)$.

\begin{lemma}
\label{lem:SLorb}
Consider the action of $\SL(s,q^g)/\GF(q)^*$ on the set of the fat $k$-subspaces of $V$.
For $k < s$, the action is transitive.
For $k = s$, $\mathcal{F}_k$ splits into $\frac{q^g - 1}{q-1}$ orbits of equal length.
\end{lemma}

\begin{Proof}
	Let $U$ be a fat $k$-subspace of $V$ and let $B$ be an ordered $\GF(q)$-basis of $U$.
	Then $B$ is an ordered $\GF(q^g)$-basis of $\langle U\rangle_{\GF(q^g)}$.

	For $k < s$, $B$ can be extended to an ordered $\GF(q^g)$-basis $B'$ of $V$.
	Let $A$ be the $(s\times s)$-matrix over $\GF(q^g)$ whose rows are given by $B'$.
	By scaling one of the vectors in $B' \setminus B$, we may assume $\det(A) = 1$.
	Now the mapping $V\to V$, $x\mapsto xA$ is in $\SL(s,q^g)$ and maps the fat 
    $k$-subspace $\langle e_1,\ldots,e_k\rangle$ to $U$ ($e_i$ denoting the $i$-th standard vector of $V$).
	Thus, the action of $\SL(s,q^g)/\GF(q)^*$ is transitive on $\mathcal{F}_k$.

	It remains to consider the case $k = s$.
	Let $A$ be the $(s\times s)$-matrix over $\GF(q^g)$ whose rows are given by $B$.
	As any two $\GF(q)$-bases of $U$ can be mapped to each other by a $\GF(q)$-linear map, we see that up to a factor in $\GF(q)^*$, $\det(A)$ does not depend on the choice of $B$.
	Thus,
	\[
		\det(U) := \det(A)\cdot\GF(q)^*\in \GF(q^g)^*/\GF(q)^*
	\]
	is invariant under the action of $\SL(s,q^g)$ on $\mathcal{F}_k$.
	It is readily checked that every value in $\GF(q^g)^*/\GF(q)^*$ appears as the invariant $\det(U)$ for some fat $s$-subspace $U$, and that two fat $s$-subspaces having the same invariant can be mapped to each other within $\SL(s,q^g)$.
	Thus, the number of orbits of the action of $\SL(s,q^g)$ on $\mathcal{F}_s$ is given by the number 
    $\#(\GF(q^g)^*/\GF(q)^*) = \frac{q^g - 1}{q-1}$ of invariants.
	As $\SL(s,q^g)$ is normal in $\GL(s,q^g)$ which acts transitively on $\mathcal{F}_s$, all orbits have the same size.
	Modding out the kernel $\GF(q)^*$ of the action yields the statement in the lemma.
\end{Proof}

\begin{theorem}
        \label{thm_qgdd_orbit_construction}
	Let $V$ be a vector space over $\GF(q)$ of dimension $gs$ with $g\geq 2$ and $s\geq 3$.
	Let $\cG$ be a Desarguesian $(g-1)$-spread in $\PG(V)$.
	\begin{enumerate}
	\item 
	For $k \in \{3,\ldots,s-1\}$, $(V,\cG,\mathcal{F}_k)$ is a $(gs,g,k,\lambda)_q$-GDD with
	\[
		\lambda = q^{(g-1)(\binom{k}{2}-1)} \prod_{i=2}^{k-1}\frac{q^{g(s-i)}-1}{q^{k-i}-1}\text{.}
	\]
	\item
	For each $\alpha\in\{1,\ldots,\frac{q^g - 1}{q-1}\}$, the union $\cB$ of any $\alpha$ orbits of the action of
    \linebreak 
    $\SL(s,q^g) / \GF(q)^*$ on $\mathcal{F}_s$ gives a $(gs,g,s,\lambda)_q$-GDD $(V,\cG,\cB)$ with
	\[
		\lambda = \alpha q^{(g-1)(\binom{s}{2}-1)}\prod_{i=2}^{s-2}\frac{q^{gi}-1}{q^i-1}\text{.}
	\]
	\end{enumerate}
\end{theorem}

\begin{Proof}
We may assume $V = \GF(q^g)^s$ and $\cG = \qbin{V}{1}{q^g}$.
The lines covered by the elements of $\cG$ are exactly the non-fat $\GF(q)$-subspaces of $V$ of dimension $2$.

Part~1:
By Lemma~\ref{lem:transitive} and Lemma~\ref{lem:SLorb}, $(V,\cG,\mathcal{F}_k)$ is a GDD.
Double counting yields $\#\mathcal{F}_2\cdot \lambda = \#\mathcal{F}_k \cdot\qbin{k}{2}{q}$.
Using Lemma~\ref{lem:fatsize}, this equation transforms into the given formula for $\lambda$.

Part~2: In the case $k = s$, by Lemma~\ref{lem:SLorb}, each union $\cB$ of $\alpha\in\{1,\ldots,\frac{q^g - 1}{q-1}\}$ orbits under the action of $\SL(s,q)/\GF(q)^*$ on $\mathcal{F}_s$ yields a GDD with
\[
\lambda
= \alpha q^{(g-1)(\binom{s}{2}-1)}\frac{q-1}{q^g - 1}  \prod_{i=2}^{s-1}\frac{q^{g(s-i)}-1}{q^{s-i}-1}
= \alpha q^{(g-1)(\binom{s}{2}-1)}\prod_{i=2}^{s-2}\frac{q^{gi}-1}{q^i-1}
\text{.}
\]
\end{Proof}

\begin{remark}
In the special case $g=2$, $k = s=3$ and $\alpha = 1$ the second case of 
Theorem~\ref{thm_qgdd_orbit_construction} yields $(6,2,3,q^2)_q$-GDDs.
These parameters match the ``residual construction for the $q$-Fano plane'' in 
\cite{EtzionHooker2017}.
\end{remark}

\begin{example}
We look at the case $g=2$, $k=s=3$ for $q=3$.
The ambient space is the $\GF(3)$-vector space $V = \GF(9)^3 \cong \GF(3)^6$.
We will use the representation $\GF(9) = \GF(3)(a)$, where $a$ is a root of the irreducible polynomial $x^2 - x - 1\in\GF(3)[x]$.

By Lemma~\ref{lem:fatsize}, out of the $\qbin{6}{3}{3} = 33880$ $3$-dimensional $\GF(3)$-subspaces of $V$,
\[
	\#\mathcal{F}_3 = 3^3 \cdot \frac{3^6 - 1}{3^3-1} \cdot \frac{3^4 -1}{3^2 - 1} \cdot 
    \frac{3^2 - 1}{3 - 1} = 27 \cdot 28 \cdot 10 \cdot 4 = 30240
\]
are fat.
According to Lemma~\ref{lem:SLorb}, the action of $\SL(3,9)/\GF(3)^*$ splits these fat subspaces $U$ 
into $4$ orbits of equal size $30240 / 4 = 7560$.
The orbits are distinguished by the invariant
\[
	\det(U) \in \GF(9)^*/\GF(3)^* = \{\{1,-1\},\;\{a,-a\},\;\{a+1,-a-1\},\;\{a-1,-a+1\}\}\text{.}
\]
The four orbits will be denoted by $O_1$, $O_a$, $O_{a+1}$ and $O_{a-1}$, accordingly.

As a concrete example, we look at the $\GF(3)$-row space $U$ of the matrix
\[
		A = \begin{pmatrix}
			a & 0 & a+1 \\
			0 & 1 & 0 \\
			0 & -a+1 & a
		\end{pmatrix}
		\in\GF(9)^{3\times 3}
\]
Then $\det(A) = a^2 = a + 1$, so $\det(U) = (a+1)\cdot\GF(3)^* = \{a+1,-a-1\}$ and thus $U \in O_{a+1}$.
Using the ordered $\GF(3)$-basis $(1,a)$ of $\GF(9)$, $\GF(9)$ may be identified with $\GF(3)^2$ and $V$ 
may be identified with $\GF(3)^6$.
The element $1\in\GF(9)$ turns into $(1,0)\in\GF(3)^2$, $a$ turns into $(0,1)$, $a-1$ turns into $(-1,1)$, etc.
The subspace $U$ turns into the row space of the matrix
	\[
		\begin{pmatrix}
			0 & 1 & 0 & 0 & 1 & 1 \\
			0 & 0 & 1 & 0 & 0 & 0 \\
			0 & 0 & 1 & -1 & 0 & 1
		\end{pmatrix}
		\in\GF(3)^{3\times 6}\text{.}
	\]
By Theorem~\ref{thm_qgdd_orbit_construction}, any disjoint union of $\alpha\in\{1,2,3,4\}$ orbits in 
$\{O_1, O_a, O_{a+1}, O_{a-1}\}$ is a $(6,2,3,9\alpha)_3$-GDD with respect to the Desarguesian line spread 
given by all $1$-dimensional $\GF(9)$-subspaces of $V$ (considered as $2$-dimensional $\GF(3)$-subspaces).
\end{example}

\begin{remark}
	A fat $k$-subspace ($k \in\{3,\ldots,s\}$) is always scattered with respect to the Desarguesian spread $\qbin{V}{1}{q^g}$.
	The converse is only true for $g = 2$.
	Thus, Theorem~\ref{thm_qgdd_orbit_construction} implies that the set of all scattered $k$-subspaces with respect to the Desarguesian line spread of $\GF(q)^{2s}$ is a $(2s,2,k,\lambda_{\max})_q$-GDD.
\end{remark}

\section{Computer constructions}\label{sec:computer}
An element $\pi \in \PGGL(v, q)$ is an automorphism of a
$(v,g,k,\lambda)_q$-GDD if
$\pi(\cG) = \cG$
and $\pi(\cB) = \cB$.

Taking the Desarguesian $(g-1)$-spread and applying the Kramer-Mesner method
\cite{KramerMesner:76} with the tools described in \cite{Braun-Kerber-Laue-2005-DCC34[1]:55-70,qdesigns2017,qdesignscomputer2017} to
the remaining blocks,
we have found $(v, g, k,\lambda)_q$-GDDs for the parameters
listed in Tables~\ref{tab:results_q2}, \ref{tab:results_q3}.
In all cases, the prescribed automorphism groups are subgroups of the 
\emph{normalizer $\langle \sigma, \phi\rangle$ of a Singer cycle group} generated by 
an element $\sigma$ of order $q^v-1$ and by the Frobenius automorphism $\phi$, 
see \cite{qdesignscomputer2017}.
Note that the presented necessary conditions for $\lambda_\delt$ turn out to be tight in several cases. 

The $q$-GDDs computed with the Kramer-Mesner approach are available in electronic form
at~\cite{BKKNW:2018}. 
The downloadable zip file contains for each parameter set $(v,k,g,q)$ a bzip2-compressed file
storing the used spread and the blocks of the $q$-GDDs for all values of $\lambda$ in the data format JSON.

\begin{table}[ht]
\caption{Existence results for $(v, g, k, \lambda)_q$-GDD for $q=2$.}\label{tab:results_q2}
\begin{center}
\begin{tabular}{rccrrll}
\hline
$v$ & $g$ & $k$ &  $\lambda_\delt$ & $\lambda_{\max}$ & $\lambda$  & comments\\  \hline
6   & 2   & 3   &  2 & 12  &  4                          & \cite{EtzionHooker2017}\\
    &     &     &    &     & 2, 4, \ldots, 12           & $\langle \sigma^7 \rangle$ \\
    &     &     &    &     & $4\alpha$, $\alpha=1, 2, 3$    & Thm. \ref{thm_qgdd_orbit_construction} \\
6   & 3   & 3   &  3 & 6   &  3, 6                       & $\langle \sigma^{21} \rangle$ \\
8   & 2   & 3   &  2 & 60  &  2, 58                      & $\langle \sigma, \phi^4\rangle$ \\
    &     &     &    &     &  4, 6, \ldots, 54, 56, 60  & $\langle \sigma, \phi\rangle$ \\
8   & 2   & 4   &  4 & 480 & 20, 40, \ldots, 480        & $\langle \sigma, \phi\rangle$ \\
    &     &     &    &     & $160\alpha$, $\alpha=1, 2, 3$    & Thm. \ref{thm_qgdd_orbit_construction} \\
8   & 4   & 3   &  7 &  42 & 7, 21, 35                   & $\langle \sigma \rangle$ \\
    &     &     &    &     & 14, 28, 42                  & $\langle \sigma, \phi \rangle$ \\
8   & 4   & 4   &  7 &  14 & 14                          & Trivial \\
9   & 3   & 3   &  1 & 118 & 2, 3, \ldots, 115, 116, 118 & $\langle \sigma, \phi\rangle$ \\
    &     &     &    &     & $16\alpha$, $\alpha=1,\ldots,16$ & Thm. \ref{thm_qgdd_orbit_construction} \\
9   & 3   & 4   & 10 &1680 & 30, 60, \ldots, 1680       & $\langle \sigma, \phi\rangle$ \\
10  & 2   & 3   & 14 & 252 & 14, 28, \ldots, 252        & $\langle \sigma, \phi\rangle$ \\
10  & 2   & 5   &  8 &     & $23040\alpha$, $\alpha=1,\ldots,3$ & Thm. \ref{thm_qgdd_orbit_construction} \\
10  & 5   & 3   & 21 & 210 & 105, 210                    & $\langle \sigma, \phi^2\rangle$ \\
12  & 2   & 3   &  2 & 1020 &  4                         & \cite{fmp:10491987} \\ 
12  & 2   & 6   & 16 &     & $12533760\alpha$, $\alpha=1,\ldots,3$ & Thm. \ref{thm_qgdd_orbit_construction} \\
12  & 3   & 4   &  2 &     & $21504\alpha$, $\alpha=1,\ldots,7$ & Thm. \ref{thm_qgdd_orbit_construction} \\
12  & 4   & 3   &  1 & 1002 & $64\alpha$, $\alpha=1,\ldots,15$ & Thm. \ref{thm_qgdd_orbit_construction} \\
\end{tabular}
\end{center}
\end{table}

\begin{table}
\caption{Existence results for $(v, g, k, \lambda)_q$-GDD for $q=3$.}\label{tab:results_q3}
\begin{center}
\begin{tabular}{rccrrll}
$v$ & $g$ & $k$ &  $\lambda_\delt$ & $\lambda_{\max}$ & $\lambda$  & comments\\  \hline
6   & 2   & 3   &  3 & 36   &  9                         & \cite{EtzionHooker2017} \\
    &     &     &    &      & $9\alpha$, $\alpha=1,\ldots,4$              & Thm. \ref{thm_qgdd_orbit_construction} \\
    &     &     &    &      & 12, 18, 24, 36             & $\langle \sigma^{13}, \phi\rangle$ \\
6   & 3   & 3   &  4 & 24   & 12, 24                     & $\langle \sigma^{14}, \phi^2\rangle$ \\
8   & 2   & 4   &  9 &  9720 & $2430\alpha$, $\alpha=1,\ldots,4$ & Thm. \ref{thm_qgdd_orbit_construction} \\
8   & 4   & 3   &  13 & 312 & 52, 104, 156, 208, 260, 312& $\langle \sigma, \phi\rangle$ \\
9   & 3   & 3   &  1  & 1077 & $81\alpha$, $\alpha=1,\ldots,13$ & Thm. \ref{thm_qgdd_orbit_construction} \\
10  & 2   & 5   &  27 & 22044960 & $5511240\alpha$, $\alpha=1,\ldots,4$ & Thm. \ref{thm_qgdd_orbit_construction} \\
12  & 2   & 6   &  81 & {\footnotesize 439267872960} & {\footnotesize $109816968240$}$\alpha$, $\alpha=1,\ldots,4$ & Thm. \ref{thm_qgdd_orbit_construction} \\
12  & 3   & 4   &  3  &     & $5373459\alpha$, $\alpha=1,\ldots,13$ & Thm. \ref{thm_qgdd_orbit_construction} \\
12  & 4   & 3   &  1 & 29472 &  $729\alpha$, $\alpha=1,\ldots,40$ & Thm. \ref{thm_qgdd_orbit_construction} \\
\end{tabular}
\end{center}
\end{table}

\begin{example}
We take the primitive polynomial $1 + x + x^3 + x^4 + x^6$, together with the 
canonical Singer cycle group generated by 

{\small
\[
\sigma = \left(\begin{array}{c}
0 1 0 0 0 0 \\
0 0 1 0 0 0 \\ 
0 0 0 1 0 0 \\ 
0 0 0 0 1 0 \\
0 0 0 0 0 1 \\
1 1 0 1 1 0
\end{array}\right)
\]}

For a compact representation we will write all $\alpha\times\beta$ matrices $X$ over $\GF(q)$ 
with entries $x_{i,j}$,
whose indices are numbered from $0$, as vectors of integers
\[
    [\sum_j x_{0,j}q^j,\ldots,\sum_j x_{\alpha-1, j}q^j], 
\]
i.e. $\sigma = [2,4,8,16,32,27]$.

The block representatives of a $(6, 2, 3, 2)_2$-GDD can be constructed by prescribing 
the subgroup $G = \langle \sigma^7\rangle$ of the Singer cycle group.
The order of $G$ is $9$, a generator is 
$
	[54,55,53,49,57,41] 
$.
The spread is generated by $[1,14]$, 
under the action of $G$ the $21$ spread elements are partitioned into $7$ orbits.
The blocks of the GDD consist of the $G$-orbits of the following $20$ generators.

{\small
\[
\begin{array}{l}
[3,16,32],
[15,16,32],
[4,8,32],
[5,8,32],
[19,24,32],
[7,24,32],
[10,4,32],\\ [1ex]
[18,28,32],
[17,20,32],
[1,28,32],
[17,10,32],
[25,2,32],
[13,6,32],
[29,30,32],\\ [1ex]
[33,12,16],
[38,40,16],
[2,36,16],
[1,36,16],
[11,12,16],
[19,20,8]
\end{array}
\]}
\end{example}

\section*{Acknowledgements}
The authors are grateful to Anton Betten who pointed out the connection to scattered subspaces and to the anonymous referee for giving helpful remarks which improved the readability paper.

This work has been performed under the auspices of
the G.N.S.A.G.A. of the C.N.R. (National Research Council)
of Italy. 

Anamari Naki\'{c} has been supported in part by the Croatian Science Foundation under the project 6732.
Deutsche Forschungsgemeinschaft (DFG) supported Sascha Kurz with KU 2430/3-1 and Alfred Wassermann
with WA1666/9-1.

\bibliographystyle{alpha} 
\newcommand{\etalchar}[1]{$^{#1}$}

\end{document}